\documentclass[11pt]{article}
\font\frak=eufm10 scaled\magstep1
\def\goth #1{\hbox{{\frak #1}}}

\def\a{\alpha}

\def\b{\beta}
\def\bw{\mbox{$\bigwedge$}}

\def\cf#1{\#(#1)}
\def\cff#1#2{\frac{#1}{\cf{#2}}}
\def\d{\delta}
\def\D{\Delta}

\def\dT{d_{\mathrm{T}}}

\def\F{\mathcal{F}}
\def\Fb{\F_{(2)}}
\def\Fba{\Fb^1}
\def\Fbb{\Fb^2}
\def\Fbc{\Fb^3}
\def\Fbd{\Fb^4}
\def\Fbe{\Fb^5}
\def\Fbf{\Fb^6}
\def\Fbg{\Fb^7}
\def\Fbj{\Fb^{10}}
\def\Fbk{\Fb^k}
\def\Fbka{\Fb^{k+1}}
\def\Fm{\F_{(m)}}
\def\Fma{\Fm^1}

\def\Fmbka{\Fm^{2k-1}}
\def\Fmk{\Fm^k}

\def\g{\gamma}

\def\half{\mbox{$\frac{1}{2}$}}

\def\l{\lambda}

\def\Maple{\textsc{Maple}}

\def\p{\partial}
\def\Pa{P_{(1)}}
\def\Pb{P_{(2)}}
\def\pd#1#2{\frac{\p #1}{\p #2}}
\def\pdb#1#2#3{\frac{\p^2 #1}{\p #2 \p #3}}

\def\Qa{Q_{(1)}}
\def\Qb{Q_{(2)}}
\def\qdot{\dot{q}}
\def\R{\mathbf{R}}

\def\smallfrac#1#2{\mbox{$\frac{#1}{#2}$}}
\def\T{\mathbf{T}}
\def\taubkE{\tau_{(2)E}^k}
\def\taubkakE{\tau_{(2)E}^{k+1,k}}

\def\ve{\varepsilon}
\def\vf#1{\frac{\p}{\p #1}}

\def\vt{\vartheta}
\def\X{\goth{X}}

\def\pd#1#2{\frac{\partial #1}{\partial #2}}
\def\pdb#1#2#3{\frac{\partial^2 #1}{\partial #2 \partial #3}}

\def\vf#1{\frac{\partial}{\partial #1}}

\def\hook{\kern 3pt \vrule height 0pt depth 0.4pt width 3pt
\vrule height 5pt depth 0.4pt\kern 3pt}

\def\bk{\rule{0em}{1.6ex}}
\def\blob{\bk\hfill\rule{0.5em}{1.5ex}}
\def\tal{\rule{0em}{2.3ex}}

\newcommand{\art}[6]{#1: #2 {\it #3\/} {\bf #4} (#5) #6}
\newcommand{\preprint}[4]{#1: #2 {\it Preprint\/} #3 (#4)}

\begin{document}

\title{The fundamental form of a second-order\\ homogeneous Lagrangian
in two variables}
\author{D. J. Saunders\\
Department of Algebra and Geometry\\
Palacky Unversity\\
779 00 Olomouc, Czech Republic\footnote{Address for
correspondence: 30 Little Horwood Road, Great Horwood,
Milton Keynes, MK17 0QE, UK}\\
e-mail david@symplectic.demon.co.uk}

\maketitle

\begin{abstract}\noindent
We construct, for a second-order homogeneous Lagrangian in two independent
variables, a differential 2-form with the property that it is closed precisely
when the Lagrangian is null. This is similar to the property of the
`fundamental Lepage equivalent' associated with first-order Lagrangians
defined on jets of sections of a fibred manifold. We show that this form may
be defined on a fourth-order frame bundle but is not, in general, projectable
to a bundle of contact elements.
\\[1ex]
{\bf Keywords:} homogeneous Lagrangian, Lepage equivalent, Euler-Lagrange
form\\[1ex]
{\bf MSC2000 Classification:} 58E99
\end{abstract}

\section{Introduction}

The `Lepage equivalents' of a Lagrangian are important tools for use when
studying variational problems on fibred manifolds:\ they are differential forms
having the same extremals as the Lagrangian form, with a further property
ensuring that their differentials give rise to the Euler-Lagrange form. If
$\pi : E \to M$ is the fibred manifold with $\dim M = m$, and if
$\l \in \Omega^m J^k\pi$ is the Lagrangian form, then any Lepage equivalent
$\theta$ of $\l$ will be defined on a jet manifold $J^l\pi$ (with, in general,
$l \geq k$) and will satisfy the conditions that $\theta - \pi_{l,k}^* \theta$
should be contact, and that for any vector field $Z \in \X(J^l\pi)$ vertical
over $E$ the contraction $i_Z d\theta$ should also be contact. The
Euler-Lagrange form $\ve$ is then the 1-contact part of $d\theta$. The basic
example of a Lepage equivalent is the Poincar\'{e} form from classical
mechanics:\ in coordinates $(t, q^a, \qdot^a)$, if $\l = L\, dt$ then
\[
\theta = L \, dt + \pd{L}{\qdot^a} (dq^a - \qdot^a dt) \, .
\]

Global Lepage equivalents may always be found for a given Lagrangian, and if
$m = 1$ then they are unique. They are never unique when $m > 1$, because
adding an arbitrary non-zero 2-contact form to any Lepage equivalent will give
a different Lepage equivalent, although such a modification will not affect the
Euler-Lagrange form. Nevertheless, when the order $k$ of the Lagrangian is no
more than 2 then it is possible to make a canonical choice of Lepage
equivalent; this cannot, however, be done when $k \geq 3$ without the
specification of some additional structure in the problem.

A particularly important question concerns the relationship between Lepage
equivalents and null Lagrangians:\ that is, Lagrangians whose Euler-Lagrange
forms vanish. Clearly if a Lagrangian has a closed Lepage equivalent then it
will be null; and when $m = 1$ then the unique Lepage equivalent of a null
Lagrangian is closed. But when $m > 1$ then a choice of Lepage equivalent
would be needed, and it is not immediately obvious how this choice should be
made.

An answer to this question for first-order Lagrangians was found by
Krupka~\cite{Krup}, and also subsequently by Betounes~\cite{Bet}. In
coordinates $x^i$ on $M$ and fibred coordinates $(x^i, u^a)$ on $E$, the Lepage
equivalent
\[
\theta = L\,\omega
+ \sum_{r=0}^{\min\{m,n\}} \frac{1}{(r!)^2}
\frac{\partial^r L}{\partial u^{a_1}_{i_1} \ldots
\partial u^{a_r}_{i_r}}
\theta^{a_1} \wedge \ldots \wedge \theta^{a_r}
\wedge \omega_{i_1 \cdots i_r},
\]
of a Lagrangian $L \omega$ (where $\omega = dx^1 \wedge \ldots \wedge dx^m$ and
$\omega_{i_1 \cdots i_r} = i_{\p / \p x^{i_r}} \omega_{i_1 \cdots i_{r-1}}$,
and where $\theta^a = du^a - u^a_i dx^i$) is closed precisely when
$L \, \omega$ is null. A similar formula for second-order (or higher-order)
Lagrangians has not yet been found, and the existence of Lepage equivalents
having this additional property has not been firmly established.

The idea of a Lepage equivalent is not directly appropriate for homogeneous
variational problems. These are problems defined on a manifold $E$ without any
given fibration over a space of independent variables, where the solution to
the variational problem is a submanifold with an orientation but without any
preferred parametrization. Instead of using jet bundles for these problems, the
Lagrangian is defined instead on the bundle of $k$-th order $m$-frames $\Fmk E$
in the manifold~\cite{CS1} (this is also called the bundle of regular $k$-th
order $m$-velocities). The Lagrangian is a function $L$ rather than an
$m$-form, and is required to satisfy a certain homogeneity condition. Factoring
the bundle of $m$-frames by the vector fields used to specify the homogeneity
condition gives rise to the bundle $J^k_+(E, m)$ of $k$-th order oriented
contact elements of dimension $m$; a Lagrangian $m$-form $\l$ on this bundle
gives rise to a homogeneous function $L$ on the frame bundle. If a fibration
$\pi : E \to M$ is given then there is an inclusion
$J^k\pi \subset J^k_+(E, m)$, and a Lagrangian form on $J^k\pi$ gives rise to
a homogeneous function $L$ on an open subset of the frame bundle.

It was shown in~\cite{CS1} that for every Lagrangian function $L$ on $\Fmk E$
it is possible to construct an $m$-form on $\Fmbka E$ called the {\em
Hilbert-Carath\'{e}odory form}\/ having the same extremals as $L$ and giving
rise to a suitable Euler-Lagrange form. The Hilbert-Carath\'{e}odory form is
projectable to the bundle of contact elements when $m = 1$ or $k \leq 2$.

It was subsequently shown in~\cite{CS2} that for a first-order Lagrangian
function there is another $m$-form on $\Fma E$ having the property that it is
closed precisely when the Lagrangian is null. This second $m$-form is
projectable to the bundle of contact elements, and if there is a fibration of
$E$ over some $m$-dimensional manifold then the restriction to the
corresponding jet bundle takes the coordinate form shown above.

The present paper is a preliminary report on a project to generalize the latter
construction to Lagrangians of arbitrary order:\ we describe a method of
constructing, for a second-order homogeneous Lagrangian in two independent
variables, a $2$-form that has the same extremals as the Lagrangian, and is
closed precisely when the Lagrangian is null. This will be the {\em fundamental
form}\/ of the Lagrangian. Although the restriction on order and dimension
suggests that this is a rather small advance, it is nevertheless significant
because the corresponding construction for second-order Lagrangians on jets of
fibrations has not yet been found.

In Section~2 of the paper we summarise the results from~\cite{CS1} that will be
needed. Section~3 contains our main theorems, and in Section~4 we investigate
the projectability of the fundamental form to a lower-order frame bundle, and
to the bundle of contact elements. Finally, Section~5 puts this work in the
context of the project as a whole, where a new calculus of vector-valued
forms~\cite{SauHom} is likely to be a useful tool in generalizing the
construction to higher orders and more variables.

\section{Homogeneous variational problems}

We consider a smooth manifold $E$ of dimension $n$, and its bundles
$\taubkE : \Fbk E \to E$ of $k$-th order $2$-frames. Important operators on
these bundles are the total derivatives and the vertical endomorphisms. The
former are vector fields $\T_i$ along the map $\taubkakE : \Fbka E \to \Fbk E$
given in coordinates by
\[
\T_i = \sum_{s=0}^k \cff{1}{i_1 \cdots i_s}
u^\a_{i i_1 \cdots i_s} \vf{u^\a_{i_1 \cdots i_s}}
\]
and the latter are type $(1,1)$ tensor fields $S^j$ on $\Fbka E$ given by
\[
S^j = \sum_{s=0}^k \cff{s+1}{i_1 \cdots i_s}
\vf{u^\a_{i i_1 \cdots i_s}} \otimes du^\a_{i_1 \cdots i_s} \, .
\]
Here and subsequently we take local coordinates $(u^\a)$ on $E$ and
corresponding jet coordinates $(u^\a_{i_1 \cdots i_s})$ on $\Fmk E$ where the
indices $i_1, \ldots, i_s$ take the values $(1, 2)$. The symbol
$\cf{i_1 \cdots i_s}$ denotes the number of distinct rearrangements of the
indices $(i_1, \ldots, i_s)$, and is needed because the jet coordinates
$(u^\a_{i_1 \cdots i_s})$ are totally symmetric in their subscripts. Intrinsic
definitions of the operators $\T_i$ and $S^j$ may be found in~\cite{CS1,
SauHom}.

We also need to use the fundamental vector fields $\D^{i_1 \cdots i_s}_i$
defined by
\[
\D^{i_1 \cdots i_s}_i = S^{i_1 \cdots i_s}(\T_i)
\]
where the tensor fields $S^i$ and $S^j$ commute, so that $S^{i_1 \cdots i_s}$
may be defined by iteration; these vector fields are well-defined on the
manifold $\Fbka E$ (rather than along the map $\taubkakE$).

We shall be interested in the actions of these objects as derivations on
differential forms. The contractions corresponding to the vector fields $\T_i$
and $\D^{i_1 \cdots i_s}_i$ will be denoted $i_i$ and $i^{i_1 \cdots i_s}_i$,
and the Lie derivatives will be denoted $d_i$ and $d^{i_1 \cdots i_s}_i$, so
that $d_i = d i_i + i_i d$ and
$d^{i_1 \cdots i_s}_i = d i^{i_1 \cdots i_s}_i + i^{i_1 \cdots i_s}_i d$. We
shall retain the symbol $S^{i_1 \cdots i_s}$ for an iterated vertical
endomorphism, and we shall also use the symbol $d_{j_1 \cdots j_s}$ for an
iterated total derivative, using the property that $d_i$ and $d_j$ commute.
Note that we often omit the pull-back maps when they would clutter up the
formul\ae.

We shall make considerable use of the commutation properties of these
operators, and so we list these properties for the action on an $r$-form; the
proofs are straightforward.

\textbf{Lemma 1}
\begin{eqnarray}
i_i d_j & = & d_j i_i \\
d_j S^i & = & S^i d_j - r \, \d^i_j \\
d^{j_1 \cdots j_s}_j S^i
& = & S^i d^{j_1 \cdots j_s}_j - \d^i_j S^{j_1 \cdots j_s} \\
i^{j_1 \cdots j_s}_j S^i
& = & S^i i^{j_1 \cdots j_s}_j + i^{i j_1 \cdots j_s}_j \\
d^k_i d_j & = & d_j d^k_i + \d^k_j d_i \\
d^{i_1 \cdots i_s}_i d_j & = & d_j d^{i_1 \cdots i_s}_i +
\sum_{r=1}^s \d^{i_r}_j d^{i_1 \cdots \widehat{i_r} \cdots i_s}_i
\end{eqnarray}
\blob

A second-order homogeneous variational problem in two independent variables is
given by a Lagrangian function $L$ on $\Fbb E$ satisfying the homogeneity
properties
\[
d^i_j L = \d^i_j L \, , \qquad d^{ik}_j L = 0 \, .
\]
Associated with such a Lagrangian are its two Hilbert forms. These are the
1-forms $\vt^i$ on $\Fbc E$ defined by
\[
\vt^i = (S^i - \half d_j S^{ji}) dL
\]
which are used to construct the Euler-Lagrange form
\[
\ve = dL - d_i \vt^i
\]
on $\Fbd E$. In coordinates
\[
\ve = \left( \pd{L}{u^\a} - d_i \left( \pd{L}{u^\a_i} \right)
+ \frac{1}{\#(ij)} d_{ij} \left( \pd{L}{u^\a_{ij}} \right)\right) du^\a \, ,
\]
incorporating the Euler-Lagrange equations for the variational problem
defined by $L$. More details of this construction may be found in~\cite{CS1}.

\section{The fundamental form}

Let $L : \Fbb E \to \R$ be a second-order homogeneous Lagrangian function. We
define the {\em fundamental form}\/ of $L$ to be the $2$-form
\[
\Theta = \Pa^2 d\vt^1 - \Pa^1 d\vt^2 \in \Omega^2 \Fbe E
\]
where the operator $\Pa^i : \Omega^2 \Fbc E \to \Omega^2 \Fbe E$ is defined by
\[
\Pa^i = \smallfrac{1}{4} S^i - \smallfrac{1}{24} d_j S^{ji}
+ \smallfrac{1}{192} d_{jk} S^{jki} \, .
\]
Our tasks in this section will be to show that $\Theta$ has the same extremals
as $L$, and that $d\Theta = 0$ precisely when the Euler-Lagrange form
$\ve = dL - d_i \vt^i$ vanishes.

We shall carry out the first task by demonstrating that a homogeneous
Lagrangian can be recovered from its fundamental form by contracting with
total derivatives. There are two stages to the argument; the first, where we
obtain the Lagrangian from the Hilbert forms, is straightforward.

\textbf{Proposition 1}
{\em
If the Lagrangian $L$ is homogeneous then
\[
L = \half i_i \vt^i \, .
\]
}

\textbf{Proof}
We use the formula for $\Pb^i$ in the definition $\vt^i = \Pb^i dL$, and the
commutation rules. Starting with
\[
i_k \vt^i = i_k (S^i - \half d_j S^{ji}) dL
\]
we have
\[
i_k S^i dL = S^i i_k dL + i^i_k dL = \d^i_k L
\]
using the homogeneity property $i^i_k dL = d^i_k L = \d^i_k L$ and the fact
that the contraction $i_k dL$ is a function so that $S^i i_k dL$ vanishes.
We also have
\[
i_k d_j S^{ji} dL = d_j i_k S^{ji} dL
= d_j (S^{ji} i_k + S^i i^j_k + S^j i^i_k + i^{ji}_k) dL = 0
\]
using the homogeneity property $i^{ji}_k dL = d^{ji}_k L = 0$ and the fact
that the other three contractions are functions and are annihilated by $S^i$.
We conclude that $i_i \vt^i = \d^i_i L =2L$.
\blob

The second stage of this argument, that we can obtain the Hilbert forms from
the fundamental form, requires considerably more work. We shall first assemble
some preliminary results.

\textbf{Lemma 2}
{\em
If the Lagrangian $L$ is homogeneous then
\[
i^i_l \vt^m = i^{ij}_l \vt^m = i^{ijk}_l \vt^m = 0 \, .
\]
}

\textbf{Proof}
We give the proof for $i^i_l \vt^m$; the remaining arguments are similar. We
have
\[
i^i_l \vt^m = i^i_l (S^m - \half d_n S^{nm}) dL \, ;
\]
and then both
\[
i^i_l S^m dL = S^m i^i_l dL + i^{im}_l dL = 0
\]
and
\[
i^i_l d_n S^{nm} dL = d_n i^i_l S^{nm} dL
= d_n (S^{nm} i^i_l + S^m i^{in}_l + S^n i^{im}_l + i^{inm}_l) dL = 0 \, ,
\]
as $i^{im}_l dL = i^{inm}_l dL = 0$ by homogeneity, and $S^m i^i_l dL = 0$
because $i^i_l dL$ is a function.

\blob

\textbf{Lemma 3}
{\em
If the Lagrangian $L$ is homogeneous then
\[
i_l \vt^m = \d^m_l L \, .
\]
}

\textbf{Proof}
We have
\[
i_l \vt^m = i_l (S^m - \half d_n S^{nm}) dL \, ;
\]
and then
\[
i_l S^m dL = S^m i_l dL + i^m_l dL = d^m_l dL = \d^m_l L
\]
using the homogeneity of $L$ and the fact that $i_l dL$ is a function so that
$S^m i_l dL$ vanishes, whereas
\[
i_l d_n S^{nm} dL = d_n i_l S^{nm} dL
= d_n (S^{nm} i_l + S^m i^n_l + S^n i^m_l + i^{nm}_l) dL = 0
\]
for similar reasons.
\blob

\textbf{Lemma 4}
{\em
If the Lagrangian $L$ is homogeneous then
\[
d^{ij}_l \vt^m = d^{ijk}_l \vt^m = 0 \, .
\]
}

\textbf{Proof}
We give the proof for $d^{ij}_l \vt^m$. We have
\[
d^{ij}_l \vt^m = d^{ij}_l (S^m - \half d_n S^{nm}) dL \, ;
\]
and then
\[
d^{ij}_l S^m dL = S^m d^{ij}_l dL - \d^m_l S^{ij} dL = - \d^m_l S^{ij} dL
\]
because $d^{ij}_l dL = 0$ by homogeneity, whereas
\begin{eqnarray*}
d^{ij}_l d_n S^{nm} dL
& = & (d_n d^{ij}_l + \d^i_n d^j_l + \d^j_n d^i_l) S^{nm} dL \\
& = & \left( d_n (S^{nm} d^{ij}_l - \d^n_n S^{ijm} - \d^m_n S^{ijn})
+ (S^{im} d^j_l - \d^i_l S^{jm} - \d^m_l S^{ji}) \right. \\
& & \quad \left. + \, (S^{jm} d^i_l - \d^j_l S^{im} - \d^m_l S^{ij}) \right) dL \\
& = & \left( (S^{im} \d^j_l - \d^i_l S^{jm} - \d^m_l S^{ji})
+ (S^{jm} \d^i_l - \d^j_l S^{im} - \d^m_l S^{ij}) \right) dL \\
& = & -2 \d^m_l S^{ij} dL
\end{eqnarray*}
where $d^i_l dL = \d^i_l dL$, $d^{ij}_l dL = 0$ by homogeneity and
$S^{ijm} dL = 0$ because $dL$ is a second-order 1-form. Consequently
$d^{ij}_l \vt^m = 0$.

The proof for $d^{ijk}_l \vt^m$ is similar, but simpler because (for instance)
$S^{ij} dL$ is replaced by $S^{ijk} dL$, so that all the terms vanish
individually.
\blob

\textbf{Lemma 5}
{\em
If the Lagrangian $L$ is homogeneous then
\[
d^i_l \vt^m = \d^i_l \vt^m - \d^m_l \vt^i \, .\\
\]
}

\textbf{Proof}
Using
\[
d^i_l \vt^m = d^i_l (S^m - \half d_j S^{jm}) dL
\]
we have
\begin{eqnarray*}
d^i_l S^m dL & = & (S^m d^i_l - \d^m_l S^i) dL \\
& = & (\d^i_l S^m - \d^m_l S^i) dL
\end{eqnarray*}
and
\begin{eqnarray*}
d^i_l d_j S^{jm} dL & = & (d_j d^i_l + \d^i_j d_l)  S^{jm} dL \\
& = & \left( d_j (S^{jm} d^i_l - \d^j_l S^{im} - \d^m_l S^{ji})
+ \d^i_j d_l  S^{jm} \right) dL \\
& = & (\d^i_l d_j S^{jm} - \d^m_l d_j S^{ji}) dL
\end{eqnarray*}
from which the result follows.
\blob

\textbf{Lemma 6}
{\em
The Hilbert forms $\vt^i$ satisfy
\[
S^i \vt^m = \half S^{im} dL
\]
so that
\[
S^i \vt^m = S^m \vt^i \, .
\]
}

\textbf{Proof}
We have
\begin{eqnarray*}
S^i \vt^m & = & S^i (S^m - \half d_j S^{jm}) dL \\
& = & \left( S^{im} - \half (d_j S^i + \d^i_j) S^{jm} \right) dL \\
& = & \half S^{im} dL
\end{eqnarray*}
because $dL$ is a second-order 1-form.
\blob

\textbf{Lemma 7}
{\em
The Hilbert forms $\vt^i$ satisfy
\[
S^i d_i \vt^m - S^m d_i \vt^i = \vt^m \, .
\]
}

\textbf{Proof}
We have
\begin{eqnarray*}
S^i d_i \vt^m - S^m d_i \vt^i
& = & (d_i S^i \vt^m + \d^i_i \vt^m) - (d_i S^m \vt^i + \d^m_i \vt^i) \\
& = & d_i (S^i \vt^m - S^m \vt^i) + \vt^m \\
& = & \vt^m \, .
\end{eqnarray*}
\blob

We now return to the relationship between the fundamental form and the Hilbert
forms.

\textbf{Proposition 2}
{\em
If the Lagrangian $L$ is homogeneous then
\[
\vt^1 = i_2 \Theta \, , \qquad \vt^2 = - i_1 \Theta \, .
\]
}

\textbf{Proof}
We use the formula for $\Pa^i$ in the definition
$\Theta = \Pa^2 d\vt^1 - \Pa^1 d\vt^2$, and the commutation rules. Starting
with
\[
i_l \Pa^i d\vt^m = i_l (\smallfrac{1}{4} S^i - \smallfrac{1}{24} d_j S^{ji}
+ \smallfrac{1}{192} d_{jk} S^{jki}) d\vt^m
\]
we again consider the terms separately. For the first term, we get
\begin{eqnarray*}
\smallfrac{1}{4} i_l S^i d\vt^m
& = & \smallfrac{1}{4} (S^i i_l + i^i_l) d\vt^m \\
& = & \smallfrac{1}{4} \left( S^i d_l \vt^m - S^i d i_l \vt^m
+ d^i_l \vt^m - d i^i_l \vt^m \right) \, ;
\end{eqnarray*}
but we know that $i^i_l \vt^m = 0$ from Lemma~2, $d i_l \vt^m = \d^m_l dL$
from Lemma~3 and $d^i_l \vt^m = \d^i_l \vt^m - \d^m_l \vt^i$ from Lemma~5,
so we get
\[
\smallfrac{1}{4} i_l S^i d\vt^m = \smallfrac{1}{4} \left( S^i d_l \vt^m
- \d^m_l S^i dL + \d^i_l \vt^m - \d^m_l \vt^i \right) \, .
\]
But we also have
\[
S^i d_i \vt^m - S^m d_i \vt^i = \vt^m
\]
from Lemma~7, so that
\begin{eqnarray*}
\smallfrac{1}{4} (i_i S^i d\vt^m - i_i S^m d\vt^i)
& = & \smallfrac{1}{4} \left( S^i d_i \vt^m
- \d^m_i S^i dL + \d^i_i \vt^m - \d^m_i \vt^i \right) \\
& & - \smallfrac{1}{4} \left( S^m d_i \vt^i
- \d^i_i S^m dL
+ \d^m_i \vt^i - \d^i_i \vt^m \right) \\
& = & \smallfrac{1}{4} \left( S^i d_i \vt^m - S^m d_i \vt^i
+ S^m dL + 2 \vt^m \right) \\
& = & \smallfrac{1}{4} \left( S^m dL
+ 3 \vt^m \right) \, .
\end{eqnarray*}

For the second term, we get
\begin{eqnarray*}
- \smallfrac{1}{24} i_l d_j S^{j i} d\vt^m
& = & - \smallfrac{1}{24} d_j i_l S^{j i} d\vt^m\\
& = & - \smallfrac{1}{24} d_j (S^j i_l + i^j_l) S^i d\vt^m \\
& = & - \smallfrac{1}{24} d_j \left( S^j (S^i i_l + i^i_l)
+ (S^i i^j_l + i^{j i}_l) \right) d\vt^m \\
& = & \smallfrac{1}{24} \left( - d_j S^{j i} d_l \vt^m
+ d_j S^{j i} d i_l \vt^m - d_j S^j d^i_l \vt^m
+ d_j S^j d i^i_l \vt^m \right) \\
& + & \smallfrac{1}{24} \left( - d_j S^i d^j_l \vt^m
+ d_j S^i d i^j_l \vt^m
- d_j d^{j i}_l \vt^m
+ d_j d i^{j i}_l \vt^m \right) \, ;
\end{eqnarray*}
but we know that $i^i_l \vt^m = 0$ from Lemma~2, $i_l \vt^m = \d^m_l L$ from
Lemma~3, $d^{ji}_l \vt^m = 0$ from Lemma~4 and
$d^i_l \vt^m = \d^i_l \vt^m - \d^m_l \vt^i$ from Lemma~5, so that
\begin{eqnarray*}
- \smallfrac{1}{24} i_l d_j S^{j i} d\vt^m
& = & \smallfrac{1}{24} \left( - d_j S^{j i} d_l \vt^m
+ d_j S^{j i} \d^m_l dL
- d_j S^j (\d^i_l \vt^m - \d^m_l \vt^i) \right. \\
& & \left. - d_j S^i (\d^j_l \vt^m - \d^m_l \vt^j) \, .
\right)
\end{eqnarray*}
Thus
\begin{eqnarray*}
& & \smallfrac{1}{24} \left( - i_i d_j S^{j i} d\vt^m
+ i_i d_j S^{j m} d\vt^i \right) \\
& = & \smallfrac{1}{24} \left( - d_j S^{j i} d_i \vt^m
+ d_j S^{j i} \d^m_i dL
- d_j S^j (\d^i_i \vt^m - \d^m_i \vt^i) \right. \\
& & \left. - d_j S^i (\d^j_i \vt^m - \d^m_i \vt^j)
\right) \\
& & - \smallfrac{1}{24} \left( - d_j S^{j m} d_i \vt^i
+ d_j S^{j m} \d^i_i dL
- d_j S^j (\d^m_i \vt^i - \d^i_i \vt^m) \right. \\
& & \left. - d_j S^m (\d^j_i \vt^i - \d^i_i \vt^j)
\right) \\
& = & \smallfrac{1}{24} \left( d_j S^{j m} d_i \vt^i
- d_j S^{j i} d_i \vt^m - d_j S^{j m} dL
- 3 d_j S^j \vt^m \right) \, ;
\end{eqnarray*}
but
\begin{eqnarray*}
d_j S^{jm} d_i \vt^i - d_j S^{ji} d_i \vt^m
& = & d_j S^j \left( S^m d_i \vt^i - S^i d_i \vt^m \right) \\
& = & - d_j S^j \vt^m
\end{eqnarray*}
from Lemma~7, so that
\[
\smallfrac{1}{24} \left( - i_i d_j S^{ji} d\vt^m
+ i_i d_j S^{j m} d\vt^i \right)
= \smallfrac{1}{24} \left( - d_j S^{jm} dL - 4 d_j S^j \vt^m \right) \, .
\]

Finally, for the third term we get
\[
\smallfrac{1}{192} i_l d_{jk} S^{jki} d\vt^m \, ,
\]
so we need to consider
\begin{eqnarray*}
i_l d_{jk} S^{jki} d\vt^m & = & d_{jk} i_l S^{jki} d\vt^m \\
& = & d_{jk} \left( S^{jki} i_l + S^{jk} i^i_l + S^{ji} i^k_l
+ S^{ki} i^j_l + S^j i^{ki}_l + S^k i^{ji}_l + S^i i^{jk}_l \right) d\vt^m
\end{eqnarray*}
where, for instance,
\[
S^{jk} i^i_l d\vt^m = S^{jk} d^i_l \vt^m - S^{jk} d i^i_l \vt^m \, .
\]
Now
$i^i_l \vt^m = i^{ki}_l \vt^m = 0$, and
\[
S^{jki} i_l d\vt^m = - \d^m_l S^{jki} dL = 0
\]
because $dL$ is a second-order 1-form. Also, $d^{ki}_l \vt^m = 0$, and we have
\[
d_{j k} S^{j k i} d_l \vt^m = 0 \\
\]
because $d_l \vt^m$ is horizontal over $\Fbb E$, and
\[
d_{j k} S^{j k} d^i_l \vt^m
= d_{j k} S^{i j} d^k_l \vt^m
= d_{j k} S^{i k} d^j_l \vt^m = 0
\]
because the $d^i_l \vt^m$ are horizontal over $\Fba E$. Thus
\[
\smallfrac{1}{192} (i_i d_{j k} S^{j k i} d\vt^m
- i_i d_{j k} S^{j k m} d\vt^i )
= 0 \, .
\]

Putting all three terms together, we now have
\[
i_i \Pa^i d\vt^m - i_i \Pa^m d\vt^i
= \smallfrac{1}{4} \left( S^m dL + 3 \vt^m \right)
+ \smallfrac{1}{24} \left( - d_j S^{j m} dL - 4 d_j S^j \vt^m \right) \, ;
\]
but
\[
S^j \vt^m = \half S^j S^m dL
\]
from Lemma~6, so we obtain
\[
- 4 d_j S^j \vt^m = - 2 d_j S^{jm} dL
\]
and consequently
\begin{eqnarray*}
i_i \Pa^i d\vt^m - i_i \Pa^m d\vt^i
& = & \smallfrac{1}{4} \left(  S^m dL
+ 3 \vt^m \right) - \smallfrac{1}{8} d_j S^{jm} dL \\
& = & \vt^m \, .
\end{eqnarray*}
We conclude, using $\Theta = \Pa^2 d\vt^1 - \Pa^1 d\vt^2$, that
$\vt^1 = i_2 \Theta$ and $\vt^2 = - i_i \Theta$.
\blob

With the help of these two propositions, we can show that a homogeneous
Lagrangian has the same extremals as its fundamental form.

\textbf{Theorem 1}
{\em
For any map $\phi : \R^2 \to E$,
\[
((j^2\phi)^* L) \, dt^1 \wedge dt^2 = (j^5\phi)^* \Theta
\]
where $t^1, t^2$ are the standard coordinates on $\R^2$; thus the two
variational problems
\[
\d \int L \, dt^1 dt^2 = 0 \, , \qquad \d \int \Theta = 0
\]
have the same extremals.
}

\textbf{Proof}
We have
\[
L = \half (i_1 \vt^1 + i_2 \vt^2) = \half (i_1 i_2 \Theta - i_2 i_1 \Theta)
= i_1 i_2 \Theta
\]
modulo pullback maps, and the result follows immediately from the properties
of contraction with total derivatives.
\blob

We now move on to our second task, to show that $d\Theta = 0$ precisely when
the Lagrangian is null. We consider the two implications separately.

\textbf{Theorem 2}
{\em
If $d\Theta = 0$ then $L$ is a null Lagrangian.
}

\textbf{Proof}
From
\[
\vt^1 = i_2 \Theta \, , \qquad \vt^2 = - i_1 \Theta \, .
\]
we have
\[
d \vt^1 = d i_2 \Theta = d_2 \Theta \, , \qquad
d \vt^2 = - d i_1 \Theta = - d_1 \Theta
\]
using $d\Theta = 0$. Then from $L = \half i_i \vt^i$ we have
\begin{eqnarray*}
dL & = & \half \left( d i_1 \vt^1 + d i_2 \vt^2 \right) \\
& = & \half \left( d_1 \vt^1 - i_1 d\vt^1 + d_2 \vt^2 - i_2 d\vt^2 \right) \\
& = & \half \left( d_1 \vt^1 + d_2 \vt^2 \right)
- \half \left( i_1 d_2 \Theta - i_2 d_1 \Theta \right) \\
& = & \half \left( d_1 \vt^1 + d_2 \vt^2 \right)
- \half \left( - d_2 \vt^2 - d_1 \vt^1 \right) \\
& = & d_1 \vt^1 + d_2 \vt^2
\end{eqnarray*}
so that $\ve = dL - d_i \vt^i = 0$ and the Lagrangian is null.
\blob

To show the converse, we must examine the relationship between the total
derivatives $d_i$ and the operators $\Pa^j$. We need the following Lemma.

\textbf{Lemma 8}
{\em
The Hilbert forms $\vt^m$ satisfy
\[
S^{ijkl} d\vt^m = 0 \, .
\]
}

\textbf{Proof}
We note first that, as $dL$ is a second-order 1-form, both $S^{ijkl} d S^m dL$
and $S^{ijk} d S^{lm} dL$ vanish; a coordinate proof of this is
straightforward. But then
\[
S^{ijkl} d d_n S^{nm} dL = S^{ijkl} d_n d S^{nm} dL
= S^{ijk} (d_n S^l + \d^l_n) S^{nm} dL = 0
\]
so that $S^{ijkl} d\vt^m = S^{ijkl} d \Pb^m dL = 0$.
\blob

We now introduce two new operators,
\[
\Qb^i : \Omega^2 \Fbd E \to \Omega^2 \Fbg E \, , \qquad
\Qa^i : \Omega^3 \Fbf E \to \Omega^3 \Fbj E \, ,
\]
by the formul\ae\
\[
\Qb^i = \smallfrac{1}{2} S^i - \smallfrac{1}{8} d_j S^{ji}
+ \smallfrac{1}{48} d_{jk} S^{jki} - \smallfrac{1}{384} d_{jkl} S^{jkli}
\]
and
\[
\Qa^i = \smallfrac{1}{6} S^i - \smallfrac{1}{54} d_j S^{ji}
+ \smallfrac{1}{648} d_{jk} S^{jki} - \smallfrac{1}{9720} d_{jkl} S^{jkli}
+ \smallfrac{1}{174960} d_{jkln} S^{jklni} \, .
\]

\textbf{Lemma 9}
{\em
The operators $\Pa^i$ and $\Qb^i$, when acting on total derivatives of the
2-form $d\vt^m$, satisfy
\begin{eqnarray}
\setcounter{equation}{1}
(\Qb^1 d_2 - d_2 \Pa^1) d\vt^m & = & 0 \\
(\Qb^2 d_1 - d_1 \Pa^2) d\vt^m & = & 0 \\
(\Qb^1 d_1 + d_2 \Pa^2) d\vt^m & = & d\vt^m \\
(\Qb^2 d_2 + d_1 \Pa^1) d\vt^m & = & d\vt^m \, .
\end{eqnarray}
}

\textbf{Proof}
As $d\vt^m$ is a 2-form, the commutation relation to use is
\[
S^i d_j - d_j S^i = 2\d^i_j \, .
\]
We prove formula (1) as an example:\ we have
\[
\Qb^1 d_2 = (\half S^1 - \smallfrac{1}{8} d_j S^{j1}
+ \smallfrac{1}{48} d_{jk} S^{jk1}
- \smallfrac{1}{384} d_{jkl} S^{jkl1}) d_2 \, ,
\]
and so
\begin{eqnarray*}
\half S^1 d_2 & = & \half d_2 S^1 \, , \\
- \smallfrac{1}{8} d_j S^{j1} d_2
& = & - \smallfrac{1}{8} d_j (d_2 S^j + 2 \d^j_2) S^1 \\
& = & - \smallfrac{1}{8} d_2 d_j S^{j1} - \smallfrac{1}{4} d_2 S^1 \, , \\
\smallfrac{1}{48} d_{jk} S^{jk1} d_2
& = & \smallfrac{1}{48} d_{jk} (d_2 S^{jk} + 2\d^j_2 S^k + 2\d^k_2 S^j) S^1 \\
& = & \smallfrac{1}{48} d_2 d_{jk} S^{jk1} + \smallfrac{1}{12} d_2 d_j S^{j1} \, , \\
- \smallfrac{1}{384} d_{jkl} S^{jkl1} d_2
& = & - \smallfrac{1}{384} d_{jkl} (d_2 S^{jkl} + 2\d^j_2 S^{kl}
+ 2\d^k_2 S^{jl} + 2\d^l_2 S^{jk}) S^1 \\
& = & - \smallfrac{1}{384} d_2 d_{jkl} S^{jkl1}
+ \smallfrac{1}{64} d_2 d_{jk} S^{jk1} \\
& = & \smallfrac{1}{64} d_2 d_{jk} S^{jk1}
\end{eqnarray*}
because $S^{jkl1} d\vt^m = 0$ by Lemma~8. Thus
\begin{eqnarray*}
\Qb^1 d_2 & = & \smallfrac{1}{4} d_2 S^1 - \smallfrac{1}{24} d_2 d_j S^{j1}
+ \smallfrac{1}{192} d_2 d_{jk} S^{jk1} \\
& = & d_2 \Pa^1 \, .
\end{eqnarray*}
The other formul\ae\ may be obtained by similar calculations.
\blob

\textbf{Lemma 10}
{\em
The operator $\Qa^i$, when acting on total derivatives of the 3-form
$d\Theta$, satisfies
\[
\Qa^i d_i d\Theta = d\Theta \, .
\]
}

\textbf{Proof}
As $d\Theta$ is a 3-form, the commutation relation to use is now
\[
S^i d_j - d_j S^i = 3\d^i_j \, .
\]
We have
\[
\Qa^i d_i = (\smallfrac{1}{6} S^i - \smallfrac{1}{54} d_j S^{ji}
+ \smallfrac{1}{648} d_{jk} S^{jki} - \smallfrac{1}{9720} d_{jkl} S^{jkli}
+ \smallfrac{1}{174960} d_{jkln} S^{jklni}) d_i \, ,
\]
and so
\begin{eqnarray*}
\smallfrac{1}{6} S^i d_i & = & \smallfrac{1}{6} (d_i S^i + 3\d^i_i) \\
& = & \smallfrac{1}{6} d_i S^i + 1 \, , \\
- \smallfrac{1}{54} d_j S^{ji} d_i
& = & - \smallfrac{1}{54} d_j (d_i S^{ji} + 3 \d^j_i S^i + 3 \d^i_i S^j) \\
& = & - \smallfrac{1}{54} d_{ji} S^{ji} - \smallfrac{1}{6} d_j S^j
\end{eqnarray*}
and so on, giving a collapsing series; the final term, involving
$S^{jklni} d\Theta$, vanishes owing to the properties of $\Theta$.
\blob

We remark that, although we have not specified in detail the properties of
$\Theta$ which result in $S^{jklni} d\Theta$ vanishing (it is similar to
Lemma~8), this doesn't really matter:\ we could instead have specified a
series with 16 terms when defining $\Qa^i$ and then, as $d\Theta$ is a 3-form
on a 5th-order frame bundle, the final term, involving
$d_i S^{j_1 \cdots j_{15} i} d\Theta$, would be guaranteed to vanish. This
affects only the omitted pull-back maps, not the final conclusion.

\textbf{Theorem 3}
{\em
If $L$ is a null Lagrangian then $d\Theta = 0$.
}

\textbf{Proof}
If $L$ is a null Lagrangian then $dL = d_i \vt^i$, and therefore
\[
\Qb^k d_i d \vt^i = \Qb^k d d_i \vt^i = \Qb^k dd L = 0 \, .
\]
But
\[
\Theta = \Pa^2 d\vt^1 - \Pa^1 d\vt^2
\]
so that
\begin{eqnarray*}
d_1 \Theta & = & d_1 \Pa^2 d\vt^1 - d_1 \Pa^1 d\vt^2 \\
& = & \Qb^2 d_1 d\vt^1 - (1 - \Qb^2 d_2) d\vt^2 \\
& = & \Qb^2 d_i d\vt^i - d\vt^2 \\
& = & -d\vt^2 \, ;
\end{eqnarray*}
similarly $d_2 \Theta = d\vt^1$. Thus
\begin{eqnarray*}
d\Theta & = & \Qa^i d_i d\Theta \\
& = & \Qa^2 dd \vt^1 - \Qa^1 dd \vt^2 \\
& = & 0.
\end{eqnarray*}
\blob

\section{Further properties of the fundamental form}

In this section we consider whether or not the fundamental form is projectable,
first to a lower-order frame bundle, and then to a bundle of oriented contact
elements.

Our first result is positive. Although the fundamental form has been defined on
a fifth-order frame bundle, it is always projectable to the fourth-order
bundle, and furthermore it is horizontal over the second-order bundle. We shall
demonstrate projectability in coordinates, using the Lie derivative action of
the locally-defined vector fields
\[
\p^{ijk}_\a = \vf{u^\a_{ijk}} \, , \qquad
\p^{ijkl}_\a = \vf{u^\a_{ijkl}} \, , \qquad
\p^{ijklm}_\a = \vf{u^\a_{ijklm}} \, ;
\]
the result will then follow from the connectedness of the fibres of
$\Fbe E \to \Fbd E$. We start with a lemma.

\textbf{Lemma 11}
{\em The Lie derivatives $\p^{ijk}_\a$, $\p^{ijkl}_\a$, $\p^{ijklm}_\a$ commute
with the contractions $S^p$. They also satisfy
\begin{eqnarray*}
\p^{ijklm}_\a d_q & = & \d^i_q \p^{jklm}_\a + \d^j_q \p^{iklm}_\a
+ \d^k_q \p^{ijlm}_\a + \d^l_q \p^{ijkm}_\a + \d^i_q \p^{ijkl}_\a \\
\p^{ijkl}_\a d_q & = & \d^i_q \p^{jkl}_\a + \d^j_q \p^{ikl}_\a
+ \d^k_q \p^{ijl}_\a + \d^l_q \p^{ijk}_\a
\end{eqnarray*}
when acting on forms or functions on $\Fbd E$ or $\Fbc E$ respectively.
}

\textbf{Proof}
The first assertion holds because, when writing a form in coordinates in terms
of the basis forms, the Lie derivatives by $\p^{ijk}_\a$, $\p^{ijkl}_\a$,
$\p^{ijklm}_\a$ affect only the coefficient functions, whereas the contractions
by $S^p$ affect only the basis forms.

The second assertion is a straightforward computation using the coordinate
expression for the total derivative $d_q$.
\blob

\textbf{Theorem 4}
{\em The fundamental form $\Theta$ is projectable to  $\Fbd E$ and is
horizontal over $\Fbb E$.
}

\textbf{Proof}
We prove the second assertion first, by showing that $S^{pqr} \Theta$ vanishes.
We have
\[
S^{pqr} S^i d\vt^m = 0
\]
by Lemma~8; for the same reason
\begin{eqnarray*}
S^{pqr} d_j S^{ji} d\vt^m
& = & \left( d_j S^{pqr} + \d^p_j S^{qr} + \d^q_j S^{pr}
+ \d^r_j S^{pq} \right) S^{ji} d\vt^m \\
& = & 0
\end{eqnarray*}
and
\begin{eqnarray*}
S^{pqr} d_{jk} S^{jki} d\vt^m
& = & \left( d_j S^{pqr} + \d^p_j S^{qr} + \d^q_j S^{pr}
+ \d^r_j S^{pq} \right) d_k S^{jki} d\vt^m \\
& = & \left( d_{jk} S^{pqr} + d_j (\d^p_k S^{qr} + \d^q_k S^{pr}
+ \d^r_k S^{pq})
+ \d^p_j (d_k S^{qr} + \d^q_k S^r + \d^r_k S^q) \right. \\
& & \quad \left. + \, \d^q_j (d_k S^{pr} + \d^p_k S^r + \d^r_k S^p)
+ \d^r_j (d_k S^{pq} + \d^p_k S^q + \d^q_k S^p)
\right) S^{jki} d\vt^m \\
& = & 0 \, ,
\end{eqnarray*}
so that $S^{pqr} \Pa^i d\vt^m = 0$ and hence $S^{pqr} \Theta = 0$. Thus $\Theta$
is horizontal over $\Fbb E$.

To show that $\Theta$ is projectable to $\Fbd E$ it is now sufficient to take a
local basis of vector fields on $\Fbe E$ vertical over $\Fbd E$ and show that
the Lie derivatives of $\Theta$ all vanish (the contractions vanish as a
consequence of the part of the theorem just proved). But, using Lemma~11, we
have
\begin{eqnarray*}
\p^{lpqrs}_\a S^i d\vt^m & = & S^i \p^{lpsrs}_\a d\vt^m \\
& = & 0
\end{eqnarray*}
because $d\vt^m$ is pulled back from $\Fbc E$, and
\begin{eqnarray*}
\p^{lpqrs}_\a d_j S^{ji} d\vt^m & = & (\d^l_j \p^{pqrs}_\a + \d^p_j \p^{lqrs}_\a
+ \d^q_j \p^{lprs}_\a + \d^r_j \p^{lpqs}_\a
+ \d^s_j \p^{lpqr}_\a) S^{ji} d\vt^m \\
& = & S^{ji} (\d^l_j \p^{pqrs}_\a + \d^p_j \p^{lqrs}_\a
+ \d^q_j \p^{lprs}_\a + \d^r_j \p^{lpqs}_\a
+ \d^s_j \p^{lpqr}_\a) d\vt^m \\
& = & 0
\end{eqnarray*}
for a similar reason. Finally,
\begin{eqnarray*}
\p^{lpqrs}_\a d_{jk} S^{jki} d\vt^m & = & (\d^l_j \p^{pqrs}_\a
+ \d^p_j \p^{lqrs}_\a + \d^q_j \p^{lprs}_\a + \d^r_j \p^{lpqs}_\a
+ \d^s_j \p^{lpqr}_\a) d_k S^{jki} d\vt^m \\
& = & \left( \d^l_j \d^p_k \p^{qrs} + \mbox{19 similar terms} \right)
S^{jki} d\vt^m \, ;
\end{eqnarray*}
but if we use the coordinate representation of the Hilbert forms $\vt^m$ as
\[
\vt^m = \left( \pd{L}{u^\g_m} - \cff{1}{mn} d_n \pd{L}{u^\g_{mn}} \right)
du^\g + \cff{1}{mn} \pd{L}{u^\g_{mn}} du^\g_n
\]
(see~\cite{CS1}) then the only non-zero terms in the expansion of
$S^{jki} d\vt^m$ are
\[
S^{jki} \left( - \cff{1}{mn} \pdb{L}{u^\b_{pl}}{u^\g_{mn}}
du^\b_{pln} \wedge du^\g
+ \cff{1}{mn} \pdb{L}{u^\b_{pl}}{u^\g_{mn}} du^\b_{pl} \wedge du^\g_n
\right) \, ,
\]
and these are second-order; thus each term of the form
$\p^{qrs}_\a S^{jki} d\vt^m$ vanishes, so that
$\p^{lpqrs}_\a d_{jk} S^{jki} d\vt^m = 0$. We conclude that
$\p^{lpqrs}_\a \Theta = 0$, so that $\Theta$ is indeed projectable to $\Fbd E$.
\blob

We shall henceforth regard $\Theta$ as being defined on $\Fbd E$, rather than
on $\Fbe E$.

Our second result is negative:\ it is {\em not} in general the case that
$\Theta$ is projectable to the manifold $J^4_+(E, 2)$ of oriented fourth-order
2-dimensional contact elements. Projectability here would require that the
contractions $i^p_l \Theta$, $i^{pq}_l \Theta$, $i^{pqr}_l \Theta$ and
$i^{pqrs}_l \Theta$, and the Lie derivatives $d^p_l \Theta$, $d^{pq}_l \Theta$,
$d^{pqr}_l \Theta$ and $d^{pqrs}_l \Theta$, should all vanish. In principle,
therefore, it would be sufficient to choose a suitable Lagrangian, substitute
into the coordinate formula for $\Theta$, and show that at least one of the
above conditions does not hold. However the coordinate formula for $\Theta$ is
already quite complicated, and it is necessary to use a homogeneous Lagrangian,
so the calculations would be rather lengthy. We shall, instead, take a more
indirect route.

\textbf{Lemma 12}
{\em The Lie derivative $d^{pqr}_s \Theta$ satisfies
\begin{eqnarray*}
d^{pqr}_s \Theta
& = & \smallfrac{5}{96} (\d^2_s S^{pqr}
+ \d^p_s S^{2qr} + \d^q_s S^{2pr} + \d^r_s S^{2pq}) d\vt^1 \\
& & \quad + \, \smallfrac{1}{96} \d^1_s (S^{2pq} d\vt^r + S^{2pr} d\vt^q
+ S^{2qr} d\vt^p) \\
& & \quad - \, \smallfrac{5}{96} (\d^1_s S^{pqr}
+ \d^p_s S^{1qr} + \d^q_s S^{1pr} + \d^r_s S^{1pq}) d\vt^2 \\
& & \quad - \, \smallfrac{1}{96} \d^2_s (S^{1pq} d\vt^r + S^{1pr} d\vt^q
+ S^{1qr} d\vt^p) \, .
\end{eqnarray*}
}

\textbf{Proof}
Once again we use the commutation relations from Lemma~1, together with some of
our other lemmas. We have
\begin{eqnarray*}
d^{pqr}_s S^i d\vt^m & = & (S^i d^{pqr}_s - \d^i_s S^{pqr}) d\vt^m \\
& = & - \d^i_s S^{pqr} d\vt^m
\end{eqnarray*}
using Lemma~4, and
\begin{eqnarray*}
d^{pqr}_s d_j S^{ij} d\vt^m & = & (d_j d^{pqr}_s + \d^p_j d^{qr}_s
+ \d^q_j d^{pr}_s + \d^r_j d^{pq}_s) S^{ij} d\vt^m \\
& = & \left( d_j (S^{ij} d^{pqr}_s - \d^i_s S^{pqrj} - \d^j_s S^{ipqr})
\right. \\
& & \quad + \, (S^{ip} d^{qr}_s - \d^i_s S^{qrp} - \d^p_s S^{iqr}) \\
& & \quad + \, (S^{iq} d^{pr}_s - \d^i_s S^{prq} - \d^q_s S^{ipr}) \\
& & \quad \left. + \, (S^{ir} d^{pq}_s - \d^i_s S^{pqr} - \d^r_s S^{ipq})
\right) d\vt^m \\
& = & - (3 \d^i_s S^{pqr} + \d^p_s S^{iqr} + \d^q_s S^{ipr} + \d^r_s S^{ipq})
d\vt^m
\end{eqnarray*}
using Lemmas~4 and~8. We also have
\begin{eqnarray*}
d^{pqr}_s d_{jk} S^{ijk} d\vt^m & = & (d_j d^{pqr}_s + \d^p_j d^{qr}_s
+ \d^q_j d^{pr}_s + \d^r_j d^{pq}_s) d_k S^{ijk} d\vt^m \\
& = & \left(\tal d_j (d_k d^{pqr}_s + \d^p_k d^{qr}_s + \d^q_k d^{pr}_s
+ \d^r_k d^{pq}_s)
\right. \\
& & \quad + \, \d^p_j (d_k d^{qr}_s + \d^q_k d^r_s + \d^r_k d^q_s)
+ \d^q_j (d_k d^{pr}_s + \d^p_k d^r_s + \d^r_k d^p_s) \\
& & \quad \left. + \, \d^r_j (d_k d^{pq}_s + \d^p_k d^q_s + \d^q_k d^p_s)
\right) S^{ijk} d\vt^m \, ;
\end{eqnarray*}
but this simplifies considerably because, for instance,
\begin{eqnarray*}
d^{pq}_s S^{ijk} d\vt^m & = & (S^{ijk} d^{pq}_s - \d^i_s S^{pqjk}
- \d^j_s S^{ipqk} - \d^k_s S^{ijpq}) d\vt^m \\
& = & 0
\end{eqnarray*}
using Lemmas~4 and~8 again, and for the same reason we have
$d^{pqr}_s S^{ijk} d\vt^m = 0$. We are left with
\begin{eqnarray*}
d^{pqr}_s d_{jk} S^{ijk} d\vt^m & = &
\left( (\d^p_j \d^q_k + \d^p_k \d^q_j) d^r_s
+ (\d^p_j \d^r_k + \d^p_k \d^r_j) d^q_s
+ (\d^q_j \d^r_k + \d^q_k \d^r_j) d^p_s \right) S^{ijk} d\vt^m \\
& = & 2 \left( S^{ipq} d^r_s - \d^i_s S^{rpq} - \d^p_s S^{irq}
- \d^q_s S^{ipr} \right. \\
& & \quad + \, S^{ipr} d^q_s - \d^i_s S^{qpr} - \d^p_s S^{iqr}
- \d^r_s S^{ipq} \\
& & \quad \left. + \, S^{iqr} d^p_s - \d^i_s S^{pqr} - \d^q_s S^{ipr}
- \d^r_s S^{iqp} \right) d\vt^m \\
& = & 2\left( S^{ipq} d^r_s + S^{ipr} d^q_s + S^{iqr} d^p_s \right. \\
& & \quad \left. - \, 3\d^i_s S^{pqr} - 2\d^p_s S^{iqr} - 2\d^q_s S^{ipr}
- 2\d^r_s S^{ipq} \right) d\vt^m \\
& = & 2\d^m_s (S^{ipq} d\vt^r + S^{ipr} d\vt^q + S^{iqr} d\vt^p) \\
& & \quad - \, 2(3\d^i_s S^{pqr} - \d^p_s S^{iqr} - \d^q_s S^{ipr}
- \d^r_s S^{ipq}) d\vt^m
\end{eqnarray*}
using Lemma~5. Putting all this together, we obtain
\begin{eqnarray*}
d^{pqr}_s \Pa^i \vt^m
& = & - \smallfrac{1}{24} \d^i_s S^{pqr} d\vt^m
+ \smallfrac{1}{24} (3 \d^i_s S^{pqr} + \d^p_s S^{iqr} + \d^q_s S^{ipr}
+ \d^r_s S^{ipq}) d\vt^m \\
& & \quad + \, \smallfrac{1}{96} \d^m_s (S^{ipq} d\vt^r + S^{ipr} d\vt^q
+ S^{iqr} d\vt^p) \\
& & \quad - \, \smallfrac{1}{96} (3\d^i_s S^{pqr} - \d^p_s S^{iqr}
- \d^q_s S^{ipr} - \d^r_s S^{ipq}) d\vt^m \\
& = & \smallfrac{5}{96} (\d^i_s S^{pqr}
+ \d^p_s S^{iqr} + \d^q_s S^{ipr} + \d^r_s S^{ipq}) d\vt^m \\
& & \quad + \, \smallfrac{1}{96} \d^m_s (S^{ipq} d\vt^r + S^{ipr} d\vt^q
+ S^{iqr} d\vt^p)
\end{eqnarray*}
from which the result follows.
\blob

\textbf{Lemma 13}
{\em If $\Theta$ is projectable to $J^4_+(E, 2)$ then
$S^{112} d\vt^1 = S^{111} d\vt^2 = 0$.
}

\textbf{Proof}
If $\Theta$ is projectable to $J^4_+(E, 2)$ then in particular we must have
$d^{111}_1 \Theta = d^{112}_2 \Theta = 0$; but from the previous Lemma we
find that
\begin{eqnarray*}
d^{111}_1 \Theta
& = & \smallfrac{3}{16} S^{112} d\vt^1 - \smallfrac{5}{24} S^{111} d\vt^2 \\
d^{112}_2 \Theta & = & \smallfrac{1}{12} S^{112} d\vt^1 -
\smallfrac{1}{16} S^{111} d\vt^2 \, .
\end{eqnarray*}
\blob

\textbf{Corollary}
{\em If $\Theta$ is projectable to $J^4_+(E, 2)$ then
\[
\pdb{L}{u^\b_{11}}{u^\a_{21}} du^\b \wedge du^\a = 0 \, .
\]
}

\textbf{Proof}
We again use the coordinate representation of the Hilbert forms $\vt^m$. The
only non-zero terms in the expansion of $S^{111} d\vt^2$ are
\[
S^{111} \left( - \cff{1}{2n} \pdb{L}{u^\b_{ij}}{u^\a_{2n}}
du^\b_{ijn} \wedge du^\a
+ \cff{1}{2n} \pdb{L}{u^\b_{ij}}{u^\a_{2n}} du^\b_{ij} \wedge du^\a_n
\right) \, ,
\]
and if this expression is to vanish then the condition of the Corollary
must hold.
\blob

Our task is now to find a Lagrangian $L$ that is homogeneous, but does not
satisfy the condition in the Corollary above. First-order homogeneous
Lagrangians are easy to find:\ for instance, any determinant
$u^\a_1 u^\b_2 - u^\a_2 u^\b_1$ is a (null) homogeneous Lagrangian. But
second-order homogeneous Lagrangians are rather more complicated, and so we
shall use the result (see~\cite{CS1}) that a Lagrangian 2-form $\l$ on
$J^2_+(E, 2)$ horizontal over $J^1_+(E,2)$ gives rise to a homogeneous
Lagrangian function $L = i_2 i_1 \rho^* \l$ on $\Fbb E$, where
$\rho : \Fbb E \to J^2_+(E, 2)$ is the projection. The fact that $\l$ is
horizontal means that the contraction $i_2 i_1$ with the total derivatives is
well-defined.

Take $E = \R^4$ with coordinates $(u^1, u^2, u^3, u^4)$, and let
\[
D^{12} = u^1_1 u^2_2 - u^1_2 u^2_1 \, , \qquad
D^{23} = u^2_1 u^3_2 - u^2_2 u^3_1 \, , \qquad
D^{34} = u^3_1 u^4_2 - u^3_2 u^4_1
\]
be the three determinants on $\Fbb E$, so that the functions
\[
F_1 = \frac{D^{23}}{D^{12}} \, , \qquad F_2 = \frac{D^{34}}{D^{12}}
\]
(defined on a suitable open submanifold) are projectable to $J^2_+(E, 2)$. We
then construct the 2-form $dF_1 \wedge dF_2$, which is also projectable to
$J^2_+(E, 2)$. The projection of this 2-form is certainly horizontal over
$J^1_+(E, 2)$, so we may define a homogeneous Lagrangian function $L$ by
\[
L = i_2 i_1 (dF_1 \wedge dF_2) \, .
\]
Of course we obtain a null Lagrangian, but this has no bearing on the argument.

We now use this Lagrangian function in a \Maple\ calculation. The coordinate
expression of $L$ involves over a page of \Maple\ output, and is of no
particular interest. But the calculations confirm, as we expect, that
$d^i_j L = \d^i_j L$ and that $d^{ij}_k L = 0$, so that $L$ is indeed
homogeneous; they also give
\[
\pdb{L}{u^1_{11}}{u^2_{12}} - \pdb{L}{u^2_{11}}{u^1_{12}}
= \frac{4 u^2_2 u^3_2 D^{34}}{(D^{12})^3}
= \frac{4 u^2_2 u^3_2 (u^3_1 u^4_2 - u^3_2 u^4_1)}{
(u^1_1 u^2_2 - u^1_2 u^2_1)^3} \not= 0 \, ,
\]
showing that the condition of the Corollary is not satisfied, and therefore
that the corresponding fundamental form $\Theta$ is not projectable to
$J^4_+(E, 2)$.

\textbf{Theorem 5}
{\em If $L$ is a homogeneous Lagrangian defined on an open submanifold of
$\Fbb E$ then its fundamental form $\Theta$, defined on the corresponding open
submanifold of $\Fbd E$, will not in general be projectable to the bundle
$J^4_+(E, 2)$ of oriented fourth-order 2-dimensional contact elements.
}
\blob

\section{Further developments}

As remarked in the Introduction, a construction for second-order Lagrangians
in two independent variables is a rather small advance; a generalization to
higher orders and more variables would be desireable. The problem, of course,
is that the calculations rapidly become unmanageable without the use of more
complicated machinery.

The extension to higher orders can be aided by the use of a simple multi-index
notation for jet variables, but including extra independent variables requires
a more sophisticated tool. This has been developed in~\cite{SauHom}, and
involves the use of certain vector-valued forms on frame bundles, namely
those taking their values in $\bw^s \R^{m*}$ for the case of $m$ independent
variables. The total derivatives can be combined into a coboundary operator
$\dT$ on the spaces of these forms, and this is (modulo pull-backs) globally
exact:\ in fact the various operators $\Pa^i$, $\Pb^i$, $\Qa^i$, $\Qb^i$ used
above are truncated components of the homotopy operator for $\dT$, and Lemmas~9
and~10 are special cases of the homotopy formula. The coboundary operator may
be combined with the exterior derivative to give, for each order, a homogeneous
variational bicomplex, and the step from the Lagrangian to the Hilbert forms
involves a diagonal move across one square of the bicomplex together with an
increase in the order. The step from the Hilbert forms to the fundamental form,
in both the general first-order case and the two-variable second-order case,
involves further diagonal moves and further increases in the order until the
edge of the appropriate bicomplex has been reached.

We can therefore see what the answer ought to be in the general higher-order
case:\ we simply carry out the procedure above. The proof that the resulting
fundamental form is closed for a null Lagrangian is then an easy consequence of
the homotopy formula. The remaining task is to demonstrate that the Lagrangian
can be reconstructed from the fundamental form by contraction with total
derivatives. This involves a computation of significant complexity, and work
on the project continues.

\end{document}